\documentclass[11pt]{article}
\usepackage{epsfig}
\newtheorem{theorem}{Theorem}
\newtheorem{lemma}{Lemma}

\begin{document}

\title{On the localization theorem in equivariant cohomology}

\author{Michel Brion \and Mich\`ele Vergne}
\date{November 7, 1997}
\maketitle

\begin{abstract}
We present a simple proof of a precise version of the
localization theorem in equivariant cohomology. As an application, we
describe the cohomology algebra of any compact symplectic
variety with a multiplicity-free action of a compact Lie group. This
applies in particular to smooth, projective spherical varieties. 
\end{abstract}

\section{A precise version of the localization theorem}

Let $X$ be a topological space with an action of a compact torus $T$.
Let $H^*_T(X)$ be the equivariant cohomology algebra of $X$ with
coefficients in the field {\bf Q} of rational numbers. The equivariant
cohomology algebra of the point is denoted by $S_T$; then $H^*_T(X)$
is a $S_T$-algebra. Any weight of $T$ defines an element of degree 2
of $S_T$; this identifies $S_T$ with the symmetric algebra over 
{\bf Q} of the group of weights of $T$.

Let $\Gamma\subset T$ be a subtorus, let $X^{\Gamma}\subset X$ be
its fixed point set and let 
$$i_{\Gamma}:X^{\Gamma}\to X,~ ~ i_{T,\Gamma}:X^T\to X^{\Gamma}$$
be the inclusion maps. They define homomorphisms of $S_T$-algebras
$$i_{\Gamma}^*:H^*_T(X)\to H^*_T(X^{\Gamma}),~ ~
i_{T,\Gamma}^*:H^*_T(X^{\Gamma})\to H^*_T(X^T).$$  
Recall that the $S_T$-algebra $H^*_T(X^{\Gamma})$ is isomorphic to
$S_T\otimes_{S_{T/\Gamma}}H^*_{T/\Gamma}(X^{\Gamma})$.
In particular, the $S_T$-module 
$H^*_T(X^T)=S_T\otimes_{{\bf Q}}H^*(X^T)$ is free.

The following statement is a variant of a result of Chang and Skjelbred 
(see \cite{4} \S 2 and also \cite{8} p. 63).

\begin{theorem}
Assume that

\noindent
(i) there exists a $T$-equivariant embedding of $X$ in the space of
a real representation of $T$, and that

\noindent
(ii) the $S_T$-module $H^*_T(X)$ is free.

\noindent
Then the homomorphism of $S_T$-algebras
$$i_T^*:H^*_T(X)\to H^*_T(X^T)$$
is injective, and its image is the intersection of the images of the
$$i_{T,\Gamma}^*:H^*_T(X^{\Gamma})\to H^*_T(X^T)$$
where $\Gamma$ ranges over all subtori of codimension one of $T$.
\end{theorem}

This reduces the description of the algebra $H^*_T(X)$ to the case
where $T$ is one-dimensional; examples will be given below. 
Chang and Skjelbred's result has been generalized to equivariant
intersection cohomology by Goresky, Kottwitz and MacPherson (see
\cite{6}); a version for equivariant Chow groups of smooth, projective
algebraic varieties is given in \cite{3}.

Assumption (i) holds if $X$ is a compact differentiable manifold. If
moreover $X$ is a Hamiltonian $T$-manifold, then (ii) holds, too (see
\cite{5}).

Our proof of Theorem 1 is based on the following consequence
of the localization theorem (see \cite{8} Theorem (III.1')).

\begin{lemma} For any subgroup $\Gamma\subset T$,
the $S_T$-algebra homomorphism
$$i_{\Gamma}^*:H^*_T(X)\to H^*_T(X^{\Gamma})$$
becomes an isomorphism after inverting finitely many weights of $T$
which restrict non trivially to $\Gamma$.
\end{lemma}

We give a proof of Lemma 1 for completeness. First we consider the
case where $X^{\Gamma}$ is empty. Then, because $X$ embeds into the
space of a representation of $T$, we can cover $X$ by finitely many
$T$-invariant subsets $X_1,\ldots,X_n$ such that every $X_j$ admits a
$T$-equivariant map $X_j\to T/\Gamma_j$ where $\Gamma_j\subset T$ is a
closed subgroup which does not contain $\Gamma$. For any $T$-invariant
subset $U_j\subset X_j$, the space $H^*_T(U_j)$ is a module over the
algebra $H^*_T(T/\Gamma_j)=H^*_{\Gamma_j}(pt)$. The latter is the
quotient of $S_T$ by its ideal generated by all weights of $T$ which
restrict trivially to $\Gamma_j$. For any $j$, choose such a weight
$\chi_j$ which restricts non trivially to $\Gamma$. Then
multiplication by $\chi_j$ is zero in $H^*_T(U_j)$. By Mayer-Vietoris,
multiplication by the product of the $\chi_j$ is zero in
$H^*_T(X)$. This implies our assertion.

In the general case, let $Y\subset X$ be a closed $T$-stable
neighborhood of $X^{\Gamma}$ in $X$. Let $Z$ be the closure of
$X\setminus Y$ in $X$; then $Z$ is $T$-stable and $Z^{\Gamma}$ is
empty. By the first step of the proof and Mayer-Vietoris, it follows
that restriction $H^*_T(X)\to H^*_T(Y)$ is an isomorphism after
inverting a finite set ${\cal F}$ of characters of $T$. Moreover,
${\cal F}$ is independent of $Y$. To conclude the proof, observe that
$H^*_T(X^{\Gamma})$ is the direct limit of the $H^*_T(Y)$.
\smallskip

\noindent
{\sl Proof of Theorem 1.} By Lemma 1 applied to $\Gamma=T$, 
the homomorphism
$$i_T^*:H^*_T(X)\to H^*_T(X^T)$$
becomes an isomorphism after inverting a finite family ${\cal F}$
of non-zero weights. We may assume that these weights are primitive
(that is, not divisible in the weight lattice). Because the
$S_T$-module $H^*_T(X)$ is free, it follows that $i_T^*$ is
injective.

As $i_T=i_{\Gamma}\circ i_{T,\Gamma}$ it is clear that the image of
$i_T^*$ is contained in the intersection of images of the
$i_{T,\Gamma}^*$. To prove the opposite inclusion, choose a basis 
$(e_k)_{k\in\kappa}$ of the free $S_T$-module $H^*_T(X)$.
For any $k\in\kappa$, let 
$$e_k^*:H^*_T(X)\to S_T$$
be the corresponding coordinate function. Then there exists a
$S_T$-linear map 
$$f_k:H^*_T(X^T)\to S_T[1/\chi]_{\chi\in{\cal F}}$$
such that $f_k\circ i_T^*=e_k^*$.

Let $\chi\in{\cal F}$; then its kernel $ker(\chi)\subset T$ is a
subtorus of codimension 1. Let $u\in H^*_T(X^{ker(\chi)})$.
By Lemma 1 applied to $\Gamma=ker(\chi)$, 
there exist a product $P_{\chi}$ of weights of $T$ which are not
multiples of $\chi$, such that $P_{\chi} u$ is in the image of
$i^*_{ker(\chi)}$. Setting $v:=i^*_{T,ker(\chi)}(u)$, it follows that
$P_{\chi} v$ is in the image of $i^*_T$.
Applying $f_k$, we obtain $P_{\chi} f_k(v)\in S_T$. Thus, the
denominator of $f_k(v)$ is not divisible by $\chi$.

If $v\in H^*_T(X^T)$ is in the intersection of the images of the
$i^*_{T,ker(\chi)}$ for all $\chi\in{\cal F}$, then
$f_k(v)\in S_T[1/\chi]_{\chi\in{\cal F}}$ but the denominator of 
$f_k(v)$ is not divisible by any element of ${\cal F}$; whence
$f_k(v)\in S_T$. It follows that
$v=i^*_T(\sum_{k\in\kappa}\,f_k(v) e_k)$ 
is in $i^*_T H^*_T(X)$.

\section{Cohomology of multiplicity-free spaces}

Let $K$ be a compact connected Lie group and let $X$ be a compact
symplectic manifold with a Hamiltonian action of $K$; let
$\Phi:X\to (Lie\,K)^*$ be the moment map. The $K$-variety $X$ is
called {\sl multiplicity-free} if the preimage under $\Phi$ of any
$K$-orbit is a unique $K$-orbit (see \cite{7} and \cite{12}
Proposition A1). Under this assumption, we will describe the
$K$-equivariant cohomology algebra of $X$. Recall that this algebra
$H^*_K(X)$ is isomorphic to  $H^*_T(X)^W$ where $T\subset K$ is a
maximal torus with Weyl group $W$. The $W$-equivariant map
$i^*_T:H^*_T(X)\to H^*_T(X^T)$
restricts to an injective homomorphism
$$\iota:H^*_K(X)\to H^*_T(X^T)^W~.$$
To apply Theorem 1, we will need the following straightforward result.

\begin{lemma} Let $X$ be a multiplicity-free $K$-variety and let
$\Gamma$ be a subtorus of $T$ with centralizer $K^{\Gamma}$ in $K$.
Then the number of connected components of $X^{\Gamma}$ is finite, and
each of them is a multiplicity-free $K^{\Gamma}$-variety.
\end{lemma}

In particular, $X^T$ is finite, and $\Phi$ induces a bijection from 
$X^T/W$ onto $\Phi(X^T)/W$. To describe the latter, choose a Weyl
chamber ${\cal C}\subset (Lie\,T)^*$ and set
$$\Phi(X^T)\cap{\cal C}=\{\lambda_1,\ldots,\lambda_n\}~.$$
For $1\leq i\leq n$, choose $x_i\in X^T$ such that
$\Phi(x_i)=\lambda_i$ and denote by $W_i$ the isotropy group of $x_i$
in $W$. Then the algebra $H^*_T(X^T)^W$ identifies with the subalgebra
of $S_T\times\cdots\times S_T$ ($n$ times) consisting of the
$(f_1,\ldots,f_n)$ such that each $f_i$ is invariant under $W_i$.

Let $X_1\subset X$ be the union of fixed point subsets of subtori of
codimension one in  $T$. Then $\Phi(X_1)\cap(Lie\,T)^*$ consists of
finitely many segments with ends in $W\{\lambda_1,\ldots,\lambda_n\}$.

\begin{theorem}
Notation being as above, the algebra
$H^*_K(X)$ is isomorphic via $\iota$ to the subalgebra of
$S_T\times\cdots\times S_T$ ($n$ times) consisting of the
$(f_1,\ldots,f_n)$ such that:

\noindent
(i) each $f_i$ is invariant under $W_i$, and 

\noindent
(ii) $f_i\equiv w(f_j)$ (mod $\lambda_i-w(\lambda_j)$)
whenever the segment $[\lambda_i,w(\lambda_j)]$ lies in
$\Phi(X_1)\cap(Lie\,T)^*$.
\end{theorem}

\noindent
{\sl Remarks.}
1) The cohomology algebra $H^*(X)$ (with coefficients
in {\bf Q}) is the quotient of $\iota H^*_K(X)$ by its ideal generated
by the $(f,f,\ldots,f)$ where $f$ is a non-zero homogeneous element of
$S_T^W$. Indeed, $H^*(X)$ is the quotient of $H^*_K(X)$ by its ideal
generated by the non-zero homogeneous elements of $H^*_K(point)$ (see
\cite{5}).
\smallskip

\noindent
2) A special $K$-equivariant cohomology class is the class $\sigma$ of
the equivariant symplectic form, and we have
$i_T^*(\sigma)=\sum_{x\in X^T}\,\Phi(x) [x]$. Thus, $\iota(\sigma)$
identifies with $(\lambda_1,\ldots,\lambda_n)$.
\smallskip

\noindent
3) Theorem 2 applies to smooth projective spherical varieties, and in
particular to complete symmetric varieties. The equivariant cohomology
algebra of the latter has been described by Bifet, DeConcini,
Littelmann and Procesi via other methods (see \cite{11} and references
therein).

\smallskip

\noindent
{\sl Proof.}
Let $\Gamma\subset T$ a subtorus of codimension one. Let $\chi$ be a
weight of $T$ with kernel $\Gamma$ and let $Y$ be a connected
component of $X^{\Gamma}$. By Lemma 2, the $K^{\Gamma}$-variety $Y$
is multiplicity-free. Set  $K^{\Gamma}/\Gamma:=L$; then $L$ is
isomorphic to $S^1$, ${\rm SU}(2)$ or ${\rm SO}(3)$. Two cases can
occur:

\noindent
1) $Y$ is two-dimensional. Then $Y$ is isomorphic to complex projective
line, and $Y^T$ consists of two fixed points $y$, $z$.
Restriction to these fixed points identifies $H^*_T(Y)$ to
the set of all $(f_y,f_z)\in S_T\times S_T$ such that
$$f_y\equiv f_z~ ~ ({\rm mod}\,\chi)~.\leqno(1)$$ 

If $\chi$ is not a root of $(K,T)$, then $L\simeq S^1$ and $\Phi(Y)$
is the segment $[\Phi(y),\Phi(z)]$. Thus, this segment lies in
$\Phi(X_1)\cap(Lie\,T)^*$. On the other hand, if $\chi$ is a root, let
$s\in W$ be the corresponding reflection and let $n_s\in K$ be a
representative of $s$. Then $Y$ is invariant under $n_s$, and 
$$H^*_T(Y)^s=\{(f_y,s(f_y))~\vert~f_y\in S_T\}$$
(indeed, $f-s(f)$ is divisible by $\chi$ for any $f\in S_T$).
\noindent

2) $Y$ is four-dimensional. Then $L$ is isomorphic to ${\rm SU}(2)$
or ${\rm SO}(3)$, and the $L$-variety $Y$ is either a
rational ruled surface, or the projectivization of a three-dimensional
complex representation of ${\rm SU}(2)$ (see \cite{9} and \cite{2}
Chapter IV, Appendix A). In the former case, $Y^T$ consists of four
points $y$, $s(y)$, $z$, $s(z)$ where $s$ is the non-trivial element
of the Weyl group of $(L,T/\Gamma)$; we may assume that the segment 
$[\Phi(y),\Phi(z)]$ lies in $\Phi(Y)\cap(Lie\,T)^*$. It is
easy to check that restriction to fixed points maps $H^*_T(Y)$ onto
the set of all $(f_y,f_{s(y)},f_z,f_{s(z)})\in S_T^4$ such that
$$f_y\equiv f_{s(y)}\equiv f_z\equiv f_{s(z)}~({\rm mod}\,\chi),~
f_y+f_{s(y)}\equiv f_z+f_{s(z)}~({\rm mod}\,\chi^2)~.\leqno(2)$$
It follows that
$$H^*_T(Y)^s=\{(f_y,f_z)\in S_T\times S_T~\vert~
f_y\equiv f_z~ ({\rm mod}\,\chi)\}.$$

In the latter case, we have similarly $Y^T=\{y,s(y),z\}$ where
$z=s(z)$, and $H^*_T(Y)$ consists of the
$(f_y,f_{s(y)}, f_z)\in S_T^3$ such that
$$f_y\equiv f_{s(y)}\equiv f_z~({\rm mod}\,\chi),~
f_y+f_{s(y)}\equiv 2f_z~({\rm mod}\,\chi^2)~.\leqno(3)$$
It follows that
$$H^*_T(Y)^s=\{(f_y,f_z)\in S_T\times S_T^s~\vert~
f_y\equiv f_z~({\rm mod}\,\chi\}~.$$
We conclude that the image of $\iota$ is defined by the congruences of
Theorem 2. Observe that the image of
$i_T^*:H^*_T(X)\to H^*_T(X^T)$ is defined by congruences of the
form (1), (2) or (3); this result is proved in \cite{3} for
equivariant Chow groups of spherical varieties.
\smallskip

\noindent
{\sl Example 1 (coadjoint orbits).} Let $X$ be the $K$-orbit of
$\lambda\in(Lie\,K)^*$; we may assume that $\lambda\in{\cal C}$. Then
$\Phi:X\to(Lie\,K)^*$ is the inclusion map, whence 
$\Phi(X^T)\cap{\cal C}=\{\lambda\}$ and 
$\Phi(X_1)\cap(Lie\,T)^*=W\cdot\lambda$. Then Theorem 2
reduces to the well-known isomorphism
$$H^*_K(K\cdot\lambda)=S_T^{W_{\lambda}}~.$$
Let $\Gamma\subset T$ be the kernel of a weight $\chi$. If $\chi$ is a
root with corresponding reflection $s_{\chi}$, then
$X^{\Gamma}=K^{\Gamma}W\cdot\lambda$ is a disjoint union of complex
projective lines containing $w\lambda$  and $s_{\chi}w\lambda$;
otherwise, $X^{\Gamma}=X^T$. By Theorem 1, we have 
$$H^*_T(K\cdot\lambda)=\{(f_{\mu})_{\mu\in W\cdot\lambda}~\vert~
f_{\mu}\in S_T,~ f_{\mu}\equiv f_{s_{\alpha}\mu}~({\rm mod}~\alpha)~
{\rm for~all~roots}~\alpha\}~.$$
Another description of this algebra is due to Arabia and Kostant-Kumar
(see \cite{1} and \cite{10}).
\smallskip

\noindent
{\sl Example 2 (complete conics).} Let $V$ be the vector space of
quadratic forms on ${\bf C}^3$, let $V^*$ be the dual space,
and let ${\bf P}={\bf P}(V)\times{\bf P}(V^*)$ be the product of their
projectivizations. Let $X\subset{\bf P}$ be the closure of the set of
classes $([A],[A^{-1}])$ where $A\in V$ is non-degenerate and
$A^{-1}\in V^*$ is the dual quadratic form. Then
$X$ is a smooth projective variety, called the space of complete
conics. Moreover, $X$ is multiplicity-free for the natural action of
the unitary group $K:=U(3)$, and $\Phi(X_1)\cap(Lie\,T)^*$ is given by
the following picture.

\begin{center}
\begin{picture}(0,0)%
\epsfig{file=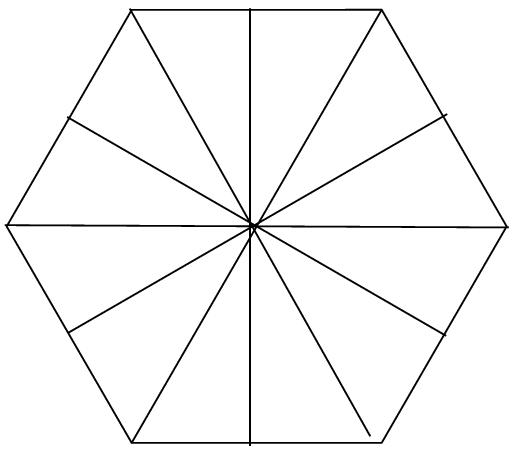}%
\end{picture}%
\setlength{\unitlength}{0.00083300in}%
\begingroup\makeatletter\ifx\SetFigFont\undefined%
\gdef\SetFigFont#1#2#3#4#5{%
  \reset@font\fontsize{#1}{#2pt}%
  \fontfamily{#3}\fontseries{#4}\fontshape{#5}%
  \selectfont}%
\fi\endgroup%
\begin{picture}(2443,2418)(2862,-2848)
\put(3976,-586){\makebox(0,0)[lb]{\smash{\SetFigFont{12}{14.4}{\rmdefault}
{\mddefault}{\updefault}$2\rho-\alpha_1$}}}
\put(4801,-736){\makebox(0,0)[lb]{\smash{\SetFigFont{12}{14.4}{\rmdefault}
{\mddefault}{\updefault}$2\rho$}}}
\put(5101,-1261){\makebox(0,0)[lb]{\smash{\SetFigFont{12}{14.4}{\rmdefault}
{\mddefault}{\updefault}$2\rho-\alpha_2$}}}
\end{picture}

\end{center}

%\begin{figure}
%\vspace{2in}
%\end{figure}

It follows that the algebra $H^*_K(X)$ consists of all triples
$(f,f_1,f_2)$ in $S_T\times S_T^{s_1}\times S_T^{s_2}$ such that
$f\equiv f_1$ (mod $\alpha_1$), $f\equiv f_2$ (mod $\alpha_2$) and
that $f_1\equiv s_{\alpha_1+\alpha_2}(f_2)$ (mod $2\alpha_1+\alpha_2$)
where $\alpha_1$, $\alpha_2$ are the simple roots, with corresponding
reflections  $s_1$, $s_2$. Other examples are given in \cite{3}.

\bigskip

\noindent
Michel BRION: Institut Fourier, B. P. 74, F-38402 Saint-Martin d'H\`eres
Cedex. e-mail {\tt mbrion@fourier.ujf-grenoble.fr}
\smallskip

\noindent
Mich\`ele VERGNE: D. M. I., \'Ecole Normale Sup\'erieure, 45 rue
d'Ulm, F-75005 Paris. e-mail {\tt vergne@dmi.ens.fr}
\end{document}